\documentclass[a4paper]{amsart}

\usepackage[T1]{fontenc}
\usepackage[latin1]{inputenc}
\usepackage{amssymb}

\usepackage[pdftitle={Comparisons between periodic, analytic and local
  cyclic cohomology},
  pdfauthor={Ralf Meyer},
  pdfsubject={Mathematics; MSC 19D55, 18G60}
]{hyperref}
\usepackage[lite]{amsrefs}

\DeclareMathOperator{\HP}{HP}
\DeclareMathOperator{\HH}{HH}
\DeclareMathOperator{\HA}{HA}
\DeclareMathOperator{\HL}{HL}
\DeclareMathOperator{\HC}{HC}
\DeclareMathOperator{\HO}{H}
\DeclareMathOperator{\K}{K}
\DeclareMathOperator{\KK}{KK}

\DeclareMathOperator{\Hom}{Hom}
\DeclareMathOperator{\Ad}{Ad}
\DeclareMathOperator{\Chern}{ch}

\newcommand*{\C}{{\mathbb{C}}}

\newcommand*{\Z}{{\mathbb{Z}}}
\newcommand*{\N}{{\mathbb{N}}}

\newcommand*{\CINF}{C^\infty}
\newcommand*{\Rat}{\mathcal{H}}
\newcommand*{\Sch}{{\mathcal{S}}}
\newcommand*{\Rep}{\mathrm{R}}

\newcommand*{\CBS}{\mathfrak{S}}

\newcommand*{\brd}{-\hspace{0pt}}
\newcommand*{\nbd}{\nobreakdash-\hspace{0pt}}

\newcommand*{\hot}{\mathbin{\hat{\otimes}}}
\newcommand*{\defeq}{\mathrel{:=}}

\theoremstyle{plain}
\newtheorem{theorem}{Theorem}[section]
\newtheorem{proposition}[theorem]{Proposition}
\newtheorem{lemma}[theorem]{Lemma}

\begin{document}

\title[Cyclic cohomology theories]{Comparisons between periodic,
  analytic and local cyclic cohomology}
\author{Ralf Meyer}
\address{Mathematisches Institut\\
         Westfälische Wilhelms-Universität Münster\\
         Einsteinstr.\ 62\\
         48149 Münster\\
         Germany
}
\email{rameyer@math.uni-muenster.de}
\urladdr{http://www.math.uni-muenster.de/u/rameyer}
\thanks{This research was supported by the EU-Network \emph{Quantum
    Spaces and Noncommutative Geometry} (Contract HPRN-CT-2002-00280)
  and the \emph{Deutsche Forschungsgemeinschaft} (SFB 478).}

\begin{abstract}
  We compute periodic, analytic and local cyclic cohomology for
  convolution algebras of compact Lie groups in order to exhibit
  differences between these theories.  A surprising result is that the
  periodic and analytic cyclic cohomology of the smooth convolution
  algebras differ, although these algebras have finite homological
  dimension.
\end{abstract}
\maketitle

\section{Introduction}
\label{sec:intro}

Cyclic cohomology is a non-commutative generalization of de Rham
homology.  Besides cyclic cohomology itself, there are several
variants like periodic cyclic cohomology, entire cyclic cohomology and
local cyclic cohomology (see \cites{Connes:Cyclic, Connes:Entire,
Puschnigg:Local}).  As in~\cite{Meyer:Analytic}, we use the more
appropriate name ``analytic cyclic cohomology'' instead of ``entire
cyclic cohomology'' in this note.  We denote these theories by
$\HP^\ast$ (periodic), $\HA^\ast$ (analytic) and $\HL^\ast$ (local),
respectively.  We compute them---together with the dual homology
theories---for some convolution algebras of compact Lie groups.  We
consider these apparently trivial examples because they clearly
exhibit differences between the three theories.

Throughout this note, $K$ is a compact Lie group.  Since we aim for
counterexamples, it would be sufficient to consider only the circle
group~$S^1$.  We allow general compact Lie groups because this creates
no additional difficulties.  Let $C^\ast(K)$ be the group
$C^\ast$\nbd{}algebra of~$K$, and let $C^\infty(K)$ and $\Rat(K)$ be the
dense subalgebras of smooth functions and of coefficients of finite
dimensional representations of~$K$, respectively.  For $K=S^1$, we get 
$$
C^\ast(S^1) \cong C_0(\Z),
\qquad
C^\infty(S^1) \cong \Sch(\Z),
\qquad
\Rat(K) \cong \bigoplus\nolimits_{n\in\Z} \C,
$$
where the product on the right hand sides is the pointwise one and
$\Sch(\Z)$ is the space of rapidly decreasing sequences.

It follows from the general properties of the local theory that the
natural maps
$$
\HL_\ast(\Rat(K)) \to
\HL_\ast(\CINF(K)) \to
\HL_\ast(C^\ast(K)) \leftarrow
\K_\ast(C^\ast(K)) \otimes_\Z \C
$$
are all isomorphisms.  The last map is the Chern-Connes character.  A
similar assertion holds in cohomology.  Thus the local theory provides
a perfect description of the de Rham (co)homology of the space of
representations~$\hat{K}$.

In order to show that the Chern-Connes character above is an
isomorphism, we prove a Universal Coefficient Theorem for the
bivariant local cyclic cohomology of $C^\ast$\nbd{}algebras.  If $A$
and~$B$ are $C^\ast$\nbd{}algebras that satisfy the UCT in bivariant
Kasparov theory, then there is a natural isomorphism
$$
\HL(A,B) \cong
\Hom (\K_\ast(A) \otimes_\Z \C, \K_\ast(B) \otimes_\Z \C)
$$
of graded complex vector spaces.

For the very small subalgebra $\Rat(K)$, all three cyclic theories
still agree, both in cohomology and homology.  This follows easily
from their behavior for direct sums, which we investigate in
Section~\ref{sec:sums}.

For $\CINF(K)$ the computation of Hochschild, cyclic and periodic
cyclic cohomology is still straightforward because $\CINF(K)$ is
biprojective.  This implies also that its homological dimension is at
most~$2$.  It is remarkable that the natural maps
$$
\HP^\ast(\Rat(K)) \leftarrow \HP^\ast(\CINF(K)),
\qquad
\HP_\ast(\Rat(K)) \to \HP_\ast(\CINF(K))
$$
fail to be isomorphisms.  Thus the periodic theory is surprisingly
sensitive to the choice of a ``smooth'' subalgebra.  If we want to
describe the de Rham homology of $\hat{K}$ by periodic cyclic
cohomology, the subalgebra $\CINF(K)$ is not yet small enough.  We
have to go down further to $\Rat(K)$.

General properties of the theories show that analytic and local cyclic
homology coincide for our three examples.  In cohomology such a
result is not available.  At the moment, direct computations of
analytic cyclic \emph{co}homology seem to be impossible.  The only
available general method to compute it is to use its homological
properties, in particular, six term exact sequences.  However, since
the periodic and analytic theory have the same homological properties,
this method can only work for algebras for which the two theories
agree.  Nevertheless, we can use the Chern-Connes character in
$\K$\nbd{}homology to construct non-trivial analytic cyclic cocycles.
We can even show that the natural map from the analytic to the local
cyclic cohomology is surjective both for $\CINF(K)$ and $C^\ast(K)$.

As a result, analytic and periodic cyclic (co)homology must differ for
$\CINF(K)$, although this algebra has finite homological dimension.
Hence Khalkhali's result that periodic and analytic cyclic cohomology
agree for a Banach algebra of finite homological dimension fails for
these very elementary nuclear Fréchet algebras.

\section{Preparations}
\label{sec:preparations}

We do not include the definitions of the cyclic theories we are
considering because that would take too much space.  We refer the reader
to \cites{Connes:Cyclic, Connes:Entire, Loday:Cyclic, Meyer:Analytic,
  Puschnigg:Local}.  All cyclic cohomology theories are defined for
algebras with additional structure like a topology or bornology.  For
the analytic and local theory, it is essential to use either
bornological algebras or inductive systems.  Since this setup may be
unfamiliar to many readers we briefly explain what bornological vector
spaces are.  See~\cite{Meyer:Analytic} for further details and
references.

A complete bornological vector space is a vector space~$V$ together
with a \emph{bornology}, which is a collection $\CBS(V)$ of subsets
that are called \emph{bounded}.  The prototypical example of a
bornology is the family of ``bounded'' subsets of a (quasi)complete
locally convex topological vector space, which we call the \emph{von
Neumann bornology}.  We equip the smooth convolution algebra
$\CINF(K)$ with this bornology.  Thus a subset $S \subseteq \CINF(K)$
is bounded if and only if the functions $D(f)\colon K\to\C$ for $f\in
S$ are uniformly bounded for any differential operator $D$ on~$K$.

Another example is the \emph{fine bornology}, which is defined on
any vector space.  It consists of the bounded subsets of finite
dimensional subspaces.  We equip the algebra $\Rat(K)$ with this
bornology.  Notice that the topology of $\Rat(K)$ is much more
complicated to describe than its bornology.

For local cyclic cohomology it is important to choose a smaller
bornology than the von Neumann bornology.  If~$A$ is a Banach algebra
or a $C^\ast$\nbd{}algebra, we \emph{always} equip it with the
\emph{precompact bornology}, which consists of precompact subsets
of~$A$.  In a complete space a subset is precompact if and only if it
is relatively compact if and only if it is contained in a compact
subset.  In particular, this is the bornology we choose on
$C^\ast(K)$.  For the Montel space $\CINF(K)$, the von Neumann
bornology and the precompact bornology coincide.  Otherwise, we should
not have chosen the von Neumann bornology on $\CINF(K)$.

We have to use the precompact bornology because approximation results
in Banach spaces usually hold uniformly on precompact subsets, but not
on von Neumann bounded subsets.  For instance, the identity map on a
Banach space with Grothendieck's approximation property can be
approximated uniformly on \emph{precompact} sets by finite rank maps.
Similarly, an asymptotic morphism is approximately multiplicative
uniformly on \emph{precompact} subsets, but usually not on von Neumann
bounded subsets.

Since we use the precompact instead of the von Neumann bornology in
the definition of analytic cyclic cohomology for Banach algebras, we
get more analytic cocycles than in the usual definition of entire
cyclic cohomology in~\cite{Connes:Entire}.  Our computations also show
that this change of bornology has a drastic effect on the resulting
theory.

The right morphisms between bornological vector spaces are the bounded
maps.  Moreover, we have an obvious notion of a bounded (bi)linear
map.  A bornological algebra is a bornological vector space with a
bounded associative multiplication.  It is straightforward to verify
that $\Rat(K)$, $\CINF(K)$ and $C^\ast(K)$ with the bornologies
defined above are complete bornological algebras.

The completed bornological tensor product~$\hot$ is defined by a
universal property for bounded bilinear maps.  We use this completed
tensor product to construct the complexes computing the various
homology theories.  For our examples we have isomorphisms of
bornological vector spaces
$$
\Rat(K) \hot \Rat(K) \cong \Rat(K\times K),
\qquad
\CINF(K) \hot \CINF(K) \cong \CINF(K\times K).
$$
The space $C^\ast(K) \hot C^\ast(K)$ is more complicated.

The periodic cyclic theory is usually defined for topological
algebras.  It is shown in~\cite{Meyer:Analytic} that for Fréchet
algebras with the precompact bornology, our bornological approach
yields the same Hochschild, cyclic and periodic cyclic (co)homology
groups as the topological approach.  We even get the same
complexes.  This equivalence between topology and bornology is based
on Grothendieck's results about compact subsets of projective tensor
products.

The bornological approach is much more flexible than the topological
approach.  There seems to be no convolution algebra on a compact Lie
group which is a Fréchet algebra and yields the correct periodic
cyclic cohomology.  We have to use $\Rat(K)$, which is best viewed as
a bornological algebra (or as an algebra without additional
structure).

A separately continuous bilinear map defined on (quasi)complete
topological vector spaces is automatically bounded.  Hence any
quasi-complete locally convex algebra with separately continuous
multiplication becomes a bornological algebra when equipped with the
von Neumann or precompact bornology.  Therefore, the bornological
tensor product is similar to the completed injective tensor product
for topological vector spaces, which is universal for separately
continuous bilinear maps.  However, it enjoys much better properties
because boundedness is a more tractable condition than separate
continuity.  For instance, the bornological tensor product is
associative, the \emph{completed} injective tensor product is only
associative under additional assumptions.

Studying bounded subsets of a bornological vector space means that we
approximate a bornological vector space from within by Banach spaces.
In fact, the basic structure theorem about complete bornological
vector spaces asserts that each such space can be written in a
canonical way as an inductive limit of Banach spaces.  Furthermore,
this construction identifies the category of complete bornological
vector spaces with a full subcategory of the category of inductive
systems of Banach spaces.  The inductive systems that come from
bornological vector spaces have the additional property that the
structure maps are all injective.

For the homological algebra behind the local cyclic theory, it is
useful to allow arbitrary inductive systems.  Therefore, Puschnigg
defines local cyclic cohomology on the category of inductive systems
of Banach algebras.  He also allows inductive systems of
\emph{admissible Fréchet algebras}.  A Fréchet algebra is called
admissible if and only if there is an open neighborhood of zero such
that for any compact subset $S \subset U$, its \emph{multiplicative
closure} $\bigcup_{n=1}^\infty S^n$ is again precompact.  For
instance, Banach algebras and smooth subalgebras of Banach algebras
(in the sense of Puschnigg~\cite{Puschnigg:Local}) are admissible
Fréchet algebras.  However, since any admissible Fréchet algebra can
be written as an inductive limit of Banach algebras, this seemingly
greater generality only makes the notation more complicated.

We remark that a complete bornological algebra need not be an
inductive limit of Banach \emph{algebras} because there may be no
multiplicatively closed bounded subsets.  For instance, the group ring
of an infinite discrete group is never an inductive limit of Banach
algebras.  It is possible, but technically difficult, to extend
Puschnigg's theory to the larger category of algebra objects in the
category of inductive systems of Banach spaces.  It is easier to
transport local cohomology to the bornological category.  This
approach yields equivalent results in all relevant applications.  We
will show this in a forthcoming article.

Fortunately, the bornological algebras that we have to consider are
all inductive limits of Banach algebras.  The Peter-Weyl theorem
asserts that $\Rat(K)$ is a countable direct sum of matrix algebras,
the sum being indexed by the irreducible representations of~$K$.  The
algebras $\CINF(K)$ and $C^\ast(K)$ are even admissible Fréchet
algebras because the multiplicative closure of any compact subset of
the open unit ball of $C^\ast(K)$ is again precompact.

We will stay away from the details of the machinery of local cyclic
cohomology, where inductive systems play a rôle.  As a result, we may
just as well work in the more intuitive setup of bornological
algebras.

\section{Behavior for direct sums}
\label{sec:sums}

Both for the computation of the theories for $\Rat(K)$ and for the
proof of the Universal Coefficient Theorem for the local cyclic
theory, we need the compatibility of the various cyclic theories with
(countable) direct sums.  The direct sum in the category of
bornological algebras is just the algebraic direct sum equipped with
the coarsest bornology making the inclusions of the summands bounded.
Hence a subset of $\bigoplus_{j\in J} A_j$ is bounded if and only if
it is contained in and bounded in a finite sum $\bigoplus_{j\in F}
A_j$ for some finite subset $F\subseteq J$.

\begin{proposition}  \label{pro:sums}
  Let $(A_j)_{j\in\N}$ be a sequence of complete bornological algebras
  and let $A \defeq \bigoplus A_j$ be the bornological direct sum.
  Let~$B$ be another bornological algebra.  Then we have:
  \begin{alignat*}{2}
    \HH_\ast(A) &\cong \bigoplus \HH_\ast(A_j),
    &\qquad
    \HH^\ast(A) &\cong \prod \HH^\ast(A_j),
    \\
    \HC_\ast(A) &\cong \bigoplus \HC_\ast(A_j),
    &\qquad
    \HC^\ast(A) &\cong \prod \HC^\ast(A_j),
    \\
    \HA_\ast(A) &\cong \bigoplus \HA_\ast(A_j),
    &\qquad
    \HA^\ast(A,B) &\cong \prod \HA^\ast(A_j,B),
    \\
    \HL_\ast(A) &\cong \bigoplus \HL_\ast(A_j),
    &\qquad
    \HL^\ast(A,B) &\cong \prod \HL^\ast(A_j,B),
  \end{alignat*}
  For the local theory, we have to assume all $A_j$ and~$B$ to be
  inductive limits of Banach algebras in order for the theory to be
  defined.
\end{proposition}

\begin{proof}
  All assertions are easy to obtain for \emph{finite} sums.  For the
  analytic and local cyclic theory, this additivity follows from
  excision.  Although Hochschild and cyclic (co)homology do not satisfy
  excision in complete generality, they do satisfy it for direct sum
  extensions (see~\cite{Loday:Cyclic}).  We will reduce countable sums
  to finite sums.
  
  The Hochschild, cyclic and analytic cyclic (co)homology of~$A$ are
  all defined as the (co)homology of certain complexes $C(A)$ associated
  to the algebra~$A$ in a natural way.  For the local theory, we have to
  take the ``local cohomology'' of the same complex as for the analytic
  theory, viewed as an inductive system of complexes.  However, since
  local cohomology enjoys even better properties than ordinary
  cohomology, this is only a notational change.  For the bivariant
  theories, we have to consider morphisms $C(A)\to C(B)$ in an
  appropriate (local) homotopy category.  Again this only creates
  notational problems.  What we really have to show is that $C(A)$ is
  chain homotopic to the direct sum of the complexes $C(A_j)$ for all
  theories we consider.  We will use that this is true for finite sums.
  The proof of additivity for finite sums actually yields this stronger
  statement in all cases.
  
  We will not define the complexes $C(A)$ here because we only use
  simple formal properties.  Besides naturality and the result for
  finite sums, we need the following compatibility with inductive
  limits: if $(A'_n)$ is a sequence of bornological algebras with
  injective structure maps $A'_n \to A'_{n+1}$, then $C(\varinjlim A'_n)
  = \varinjlim C(A'_n)$.
  
  We apply this to $A'_n \defeq \bigoplus_{j=1}^n A_j$.  Since~$A$ is
  the bornological direct limit of $(A'_n)$, we conclude that $C(A)$
  is the bornological direct limit of the complexes $C(A'_n)$ for
  $n\to\infty$.  Since~$A'_{n-1}$ is a retract of~$A'_n$, the complex
  $C(A'_{n-1})$ is a direct summand of $C(A'_n)$.  Write $C(A'_n) =
  C(A'_{n-1}) \oplus X_n$.  Then we have $C(A'_n) = X_1 \oplus \cdots
  \oplus X_n$.  As a result, $C(A) \cong \bigoplus X_n$.  Additivity
  for finite sums implies that the complexes $C(A_n)$ and~$X_n$ are
  chain homotopic.  Thus $C(A)$ is chain homotopic to the direct sum
  $\bigoplus C(A_j)$.  The assertions of the proposition follow.
\end{proof}

We do not have a general assertion about direct sums in the second
variable of the analytic and local cyclic theories.  This should be
expected because direct sums and products do not commute.

It is remarkable that the periodic theory may behave quite badly with
respect to direct sums.  The periodic cyclic cohomology is isomorphic
to the inductive limit of the cyclic cohomology groups with respect to
the periodicity operator~$S$.  Since the cyclic cohomology of a direct
sum is a direct product, the periodic cyclic cohomology of a direct
sum is an inductive limit of direct products.  However, direct
products and inductive limits do not commute!

\section{The Universal Coefficient Theorem}
\label{sec:UCT}

The local cyclic theory has very good functorial properties on the
category of $C^\ast$\nbd{}algebras.  This is the source of the Universal
Coefficient Theorem (UCT) below.  It suffices to state it for the
bivariant theory because
$$
\HL_\ast(A) \defeq \HL_\ast(\C,A),
\qquad
\HL^\ast(A) \defeq \HL_\ast(A,\C).
$$
Recall that $C^\ast$\nbd{}algebras are always equipped with the
precompact bornology.

\begin{theorem}[Universal Coefficient Theorem]
  \label{the:UCT_HL}
  Let $A$ and~$B$ be separable $C^\ast$\nbd{}algebras that satisfy the
  Universal Coefficient Theorem in Kasparov's $\KK$\brd{}theory.  Thus
  $A$ and~$B$ are $\KK$\brd{}equivalent to commutative
  $C^\ast$\nbd{}algebras.

  Then there are natural isomorphisms
  \begin{multline*}
    \HL(A,B)
    \cong
    \Hom\bigl( \HL_\ast(A), \HL_\ast(B)\bigr)
    \cong
    \Hom\bigl( \K_\ast(A)\otimes_\Z\C,
    \K_\ast(B)\otimes_\Z\C\bigr).
  \end{multline*}
\end{theorem}

\begin{proof}
  First we prove that the Chern-Connes character is an isomorphism
  \begin{equation}  \label{eq:HL_UCT_homology}
    \K_\ast(B) \otimes_\Z \C
    \cong
    \KK_\ast(\C,B) \otimes_\Z \C
    \cong
    \HL_\ast(\C,B)
    =
    \HL_\ast(B)
  \end{equation}
  provided that~$B$ satisfies the UCT in $\KK$\brd{}theory.  This
  already proves that the second map in the theorem is an isomorphism.
  There is a bivariant Chern-Connes character $\KK\to\HL$
  because~$\HL$ is a stable exact homotopy functor in both variables
  and $\KK$\brd{}theory is universal among such theories.  Invariance
  for continuous homotopies and $C^\ast$\nbd{}stability for~$\HL$ are
  proved in~\cite{Puschnigg:Local}.  Excision is proved
  in~\cite{Puschnigg:Excision} (see also \cites{Meyer:Excision,
  Meyer:Analytic}).
  
  Since a separable $C^\ast$\nbd{}algebra satisfies the UCT if and
  only if it is $\KK$\brd{}equivalent to a commutative
  $C^\ast$\nbd{}algebra, it suffices to prove the assertion for~$B$ in
  the \emph{bootstrap category} (see~\cite{Blackadar:KK}).  This is
  the smallest class of separable \emph{nuclear}
  $C^\ast$\nbd{}algebras that
  \begin{itemize}
    \item contains~$\C$;
    \item is closed under $\KK$\brd{}equivalence;
    \item is closed under countable direct sums;
    \item is closed under extensions in the following sense:
       if two $C^\ast$\nbd{}algebras in an extension belong to the
       class, so does the third.
  \end{itemize}
  In~\cite{Blackadar:KK}, the third condition is replaced by the
  seemingly stronger requirement of being closed under countable
  inductive limits.  However, a mapping telescope argument shows that
  inductive limits can be reduced to direct sums and pull backs and
  pull backs can be reduced to extensions
  (see~\cite{Schochet:Axiomatic}).
  
  Hence it suffices to show that the class of separable nuclear
  $C^\ast$\nbd{}algebras for which~\eqref{eq:HL_UCT_homology} holds has
  these properties.  Since $\HL_0(\C,\C)=\C$ and $\HL_1(\C,\C)=0$, it
  certainly contains~$\C$.  It is also closed under $\KK$\brd{}equivalence
  by the existence of the bivariant Chern-Connes character.  
  It is closed under countable direct sums because
  $$
  \HL(\C, \bigoplus\nolimits^{C^\ast} B_j) \cong
  \bigoplus \HL(\C,B_j)
  $$
  for any sequence of $C^\ast$\nbd{}algebras~$(B_j)$.  Here
  $\bigoplus^{C^\ast} B_j$ denotes the direct sum in the category of
  $C^\ast$\nbd{}algebras.  This direct sum is $\HL$\brd{}equivalent to the
  bornological direct sum by Theorem~3.15 of~\cite{Puschnigg:Local}.
  Since the finite sums are retracts of the infinite sum, the technical
  assumption for that theorem is satisfied.  Hence the claim follows
  from Proposition~\ref{pro:sums}.
  
  An extension of separable nuclear $C^\ast$\nbd{}algebras automatically
  has a bounded, even completely positive, linear section.  Hence it
  gives rise to six term exact sequences both in $\K_\ast \otimes_\Z \C$
  and~$\HL_\ast$.  A diagram chase shows that our class is closed under
  extensions.  As a result, \eqref{eq:HL_UCT_homology} holds for all~$B$
  in the bootstrap category.
  
  Next we consider the class of separable nuclear
  $C^\ast$\nbd{}algebras~$A$ for which the product in local cyclic
  cohomology gives rise to an isomorphism
  $$
  \HL(A,B) \to \Hom (\HL_\ast(A), \HL_\ast(B))
  $$
  for all bornological algebras~$B$.  We claim that this class
  contains the bootstrap category.  The proof of this claim will
  complete the proof of the theorem.
  
  The computation of $\HL_\ast(\C)$ shows that~$\C$ belongs to this
  class.  Again the class is closed under $\KK$\brd{}equivalence because
  of the existence of the bivariant Chern-Connes character.  The same
  arguments as above treat direct sums and extensions.  Hence the class
  contains the bootstrap category.
\end{proof}

We remark that if the UCT is true, then the Chern-Connes character
$$
\KK(A,B) \otimes_\Z \C \to \HL(A,B)
$$
is always injective and has dense range in an appropriate sense.
Nevertheless, it fails to be an isomorphism for $A=C^\ast(K)$ and
$B=\C$.

\section{Computation of local cyclic cohomology}
\label{sec:local}

Since $C^\ast(K)$ is a $C^\ast$\nbd{}direct sum of matrix algebras, it
belongs to the bootstrap category.  Hence the UCT of the previous
section applies.  The $\K$\nbd{}theory of $C^\ast(K)$ vanishes in odd
degrees and is isomorphic to the representation ring $\Rep(K)$ of~$K$
in even degrees.  Hence
\begin{alignat*}{2}
  \HL_0(C^\ast(K)) &\cong \Rep(K) \otimes_\Z \C,
  &\qquad
  \HL_1(C^\ast(K)) &\cong \{0\},
  \\
  \HL^0(C^\ast(K)) &\cong \Hom_\Z (\Rep(K),\C),
  &\qquad
  \HL^1(C^\ast(K)) &\cong \{0\}.
\end{alignat*}
Since~$\Rep(K)$ has a countable basis, we can identify
$\HL_0(C^\ast(K))$ and $\HL^0(C^\ast(K))$ with the direct sum and the
direct product of countably many copies of~$\C$.

\begin{proposition}  \label{pro:HL_compute}
  The inclusions $\Rat(K) \subseteq \CINF(K) \subseteq C^\ast(K)$ are
  $\HL$\brd{}equivalences and hence induce isomorphisms in local cyclic
  (co)homology.
\end{proposition}

\begin{proof}
  We claim that Theorem~3.2 in~\cite{Puschnigg:Local} applies to our
  situation.  It yields that the inclusions $\Rat(K)\to C^\ast(K)$ and
  $\Rat(K) \to \CINF(K)$ are isomorphisms in his ``stable strongly
  almost multiplicative homotopy category.''  Hence they are
  $\HL$\brd{}equivalences and the proposition follows.
  
  To apply Puschnigg's result, observe that $\Rat(K)$ is a direct sum
  of finite dimensional matrix algebras and hence the increasing union
  of a sequence $(V_n)$ of finite dimensional sub\emph{algebras}.  We
  have already remarked above that $C^\ast(K)$ and $\CINF(K)$ are
  separable admissible Fréchet algebras.  We also need to know that
  they possess Grothendieck's approximation property.  This is clear
  because $C^\ast(K)$ is nuclear as a $C^\ast$\nbd{}algebra and
  $\CINF(K)$ is nuclear as a locally convex vector space.
\end{proof}

Thus the computation of the local theory for our three algebras
follows immediately from the general properties of the theory.

Since we can compute the local cyclic (co)homology of~$\Rat(K)$ directly
using Proposition~\ref{pro:sums}, we can also use
Proposition~\ref{pro:HL_compute} backwards to compute the local cyclic
(co)homology of $C^\ast(K)$ without appealing to the Universal
Coefficient Theorem.

\section{Computation of periodic cyclic cohomology}
\label{sec:periodic}

The Peter-Weyl Theorem yields that $\Rat(K)$ is a direct sum of matrix
algebras, one for each representation of~$K$.  Using Morita invariance
and Proposition~\ref{pro:sums}, we find that its Hochschild (co)homology
vanishes except in dimensions~$\ge1$ and
$$
\HH_0(\Rat(K)) \cong \Rep(K) \otimes_\Z \C,
\qquad
\HH^0(\Rat(K)) \cong \Hom_\Z(\Rep(K), \C).
$$
Thus Connes's SBI-sequence, which relates Hochschild and cyclic
(co)homology, is highly degenerate and yields isomorphisms
\begin{alignat*}{2}
  \HC_{2n}(\Rat(K)) &\cong \Rep(K) \otimes_\Z \C,
  &\qquad
  \HC_{2n+1}(\Rat(K)) &\cong \{0\},
  \\
  \HC^{2n}(\Rat(K)) &\cong \Hom_\Z (\Rep(K),\C),
  &\qquad
  \HC^{2n+1}(\Rat(K)) &\cong \{0\},
\end{alignat*}
for all $n\in\N$, the periodicity operator~$S$ being the identity.
Hence
\begin{alignat*}{2}
  \HP_0(\Rat(K)) &\cong \Rep(K) \otimes_\Z \C,
  &\qquad
  \HP_1(\Rat(K)) &\cong \{0\},
  \\
  \HP^0(\Rat(K)) &\cong \Hom_\Z (\Rep(K),\C),
  &\qquad
  \HP^1(\Rat(K)) &\cong \{0\}.
\end{alignat*}
As a result,
$$
\HP_\ast (\Rat(K)) \cong \HL_\ast(\Rat(K)) \cong
\K_\ast(C^\ast(K)) \otimes_\Z \C
$$
and dually for cohomology.

The cyclic type homology theories of $\CINF(K)$ are similarly easy to
compute.  This computation is also a trivial special case of results
of Victor Nistor and others on the cyclic homology of crossed
products.  The following lemma asserts that $\CINF(K)$ is a
biprojective Fréchet algebra (see \cites{Baehr:Stability,
Helemskii:Homology}).

\begin{lemma}  \label{lem:CINF_biprojective}
  The map
  \begin{multline*}
    \sigma \colon
    \CINF(K) \to \CINF(K) \hot \CINF(K) \cong \CINF(K \times K),
    \qquad
    \sigma f(g_0,g_1) \defeq f(g_0 \cdot g_1),
  \end{multline*}
  is a $\CINF(K)$\brd{}bimodule section for the multiplication map
  $$
  m \colon \CINF(K) \hot \CINF(K) \to \CINF(K),
  \qquad
  f_1 \otimes f_2 \mapsto f_1 \ast f_2.
  $$
  Thus $\CINF(K)$ is projective as a $\CINF(K)$\brd{}bimodule.
  
  Its commutator quotient is isomorphic to the space $\CINF(K/\Ad K)$
  of smooth functions on~$K$ that are constant on conjugacy classes,
  that is, invariant with respect to the adjoint action.
\end{lemma}

\begin{proof}
  See~\cite{Meyer:Analytic} for the isomorphism $\CINF(K) \hot \CINF(K)
  \cong \CINF(K \times K)$.  It is clear that~$\sigma$ is a well-defined
  bounded linear map.  It is a section for~$m$ because the convolution
  takes the form $m(h)(g) = \int_K h(g g_1,g_1^{-1})\,dg_1$ on $\CINF(K
  \times K)$.  The bimodule structure on $\CINF(K)$ is defined by left
  and right convolution, the bimodule structure on $\CINF(K) \hot
  \CINF(K)$ is defined by left convolution on the left and right
  convolution on the right factor.  These actions of $\CINF(K)$ are
  integrated forms of group actions of~$K$ on these spaces, defined by
  $\lambda_g f(g_1) = f(g^{-1} g_1)$, $\lambda_g h(g_0,g_1) \defeq
  h(g^{-1} g_0,g_1)$, $\rho_g f(g_1) = f(g_1 g)$ and $\rho_g h(g_0,g_1)
  \defeq h(g_0, g_1 g)$.  It is trivial to check that~$\sigma$ is
  equivariant with respect to these group actions.  Hence it is a
  bimodule map with respect to their integrated forms.
  
  To prove that $\CINF(K)$ is projective as a bimodule, we
  iterate~$\sigma$ and consider the corresponding map
  $$
  \CINF(K)
  \to
  \CINF(K^3)
  \subseteq
  \CINF(K)^+ \hot \CINF(K) \hot \CINF(K)^+.
  $$
  Here $\CINF(K)^+$ is the space obtained by adjoining a unit to
  the non-unital algebra $\CINF(K)$.  This map is still a bimodule map
  and a section for the multiplication map.  Since $\CINF(K)^+ \hot
  \CINF(K) \hot \CINF(K)^+$ is a prototypical free bimodule,
  $\CINF(K)$ is projective.

  To compute the commutator quotient of $\CINF(K)$, we consider
  $\CINF(K)$ as a direct summand of $\CINF(K)^+ \hot \CINF(K)$ via the
  map $\sigma \circ m$.  The commutator quotient of the latter
  bimodule is canonically isomorphic to $\CINF(K)$.  A computation
  shows that the map $\sigma \circ m$ induces on commutator quotients
  the map
  $$
  P \colon \CINF(K) \to \CINF(K),
  \qquad Pf(g) \defeq \int_K f(h^{-1}gh) \,dh.
  $$
  This is exactly the projection onto the space of invariants under
  the adjoint action.  Clearly, a smooth function is invariant under
  the adjoint action if and only if it is constant on conjugacy
  classes.
\end{proof}

The lemma implies that $\CINF(K)$ is $H$\nbd{}unital (see
also~\cite{Block-Getzler-Jones:Cyclic-II}).  Thus the Hochschild
homology of $\CINF(K)$ can be computed from a projective bimodule
resolution of $\CINF(K)$, we do not have to use $\CINF(K)^+$.  Since
$\CINF(K)$ is itself a projective bimodule, there is a trivial
projective resolution.  Consequently, Hochschild homology and cohomology
vanish in dimensions $\ge 1$ and
$$
\HH_0(\CINF(K)) \cong \CINF(K/\Ad K),
\qquad
\HH^0(\CINF(K)) \cong \CINF(K/\Ad K)'.
$$
We proceed as for $\Rat(K)$ to compute the cyclic and periodic
cyclic theories and get
\begin{alignat*}{2}
  \HP_0(\CINF(K)) &\cong \CINF(K/\Ad K),
  &\qquad
  \HP_1(\CINF(K)) &\cong \{0\},
  \\
  \HP^0(\CINF(K)) &\cong \CINF(K/\Ad K)',
  &\qquad
  \HP^1(\CINF(K)) &\cong \{0\}.
\end{alignat*}
Thus the periodic and local cyclic (co)homology of $\CINF(K)$ differ
and the image of $K_\ast(C^\ast(K))\otimes_\Z \C$ under the
Chern-Connes character is a proper dense subspace of
$\HP_\ast(\CINF(K))$.

One can show that the canonical map
$$
\HP^0(\CINF(K))
\to
\HL^0(\CINF(K))
\cong
{\textstyle \prod_{n\in\N} \C}
$$
is a bijection onto the space of sequences of at most polynomial
growth.  Therefore, $\CINF(K/\Ad K) \cong \Sch$ as a bornological vector
space.  This follows also from the structure theory of nuclear
topological vector spaces because $\CINF(K/\Ad K)$ is both a subspace
and a quotient of $\CINF(K)$.

To estimate the homological dimension of the algebra $\CINF(K)$, we
have to consider resolutions of $\CINF(K)^+$, not of $\CINF(K)$.
Nevertheless, biprojectivity implies that the homological dimension is
at most~$2$ (see also~\cite{Helemskii:Homology}).

\begin{lemma}  \label{lem:CINF_dimension}
  The $\CINF(K)$\brd{}bimodule $\CINF(K)^+$ has a projective bimodule
  resolution of length~$2$.  Thus the homological dimension of the
  algebra $\CINF(K)$ is at most~$2$.
\end{lemma}

\begin{proof}
  Since $\CINF(K)^+$ is an extension of the trivial bimodule~$\C$ by
  the projective bimodule $\CINF(K)$, it suffices to construct a
  resolution of length~$2$ for~$\C$.  This resolution can then be
  patched together with the trivial projective resolution of
  length~$0$ of $\CINF(K)$.  Abbreviate $A\defeq \CINF(K)$.  The
  resolution of~$\C$ is simply the tensor product of two copies of the
  extension $A\to A^+ \to \C$, that is,
  $$
  \C
  \longleftarrow
  A^+ \hot A^+
  \longleftarrow
  A \hot A^+ \oplus A^+ \hot A
  \longleftarrow
  A \hot A
  \longleftarrow
  0.
  $$
  The bimodule structure on each summand is the obvious one.
  Since~$A$ is projective as a bimodule, it is \emph{a fortiori}
  projective as a left or right module.  Hence this is a projective
  resolution of the trivial bimodule of length~$2$.
\end{proof}

We do not really care about the periodic cyclic (co)homology of
$C^\ast(K)$ because we do not expect the periodic theory to yield good
results for $C^\ast$\nbd{}algebras, anyway.  Nevertheless, we give the
result in cohomology because we will use it below and because it is
easy to obtain using the amenability of $C^\ast(K)$.  Moreover, since
the space of irreducible representations of $C^\ast(K)$ is
$0$\nbd{}dimensional, the result actually is not too bad in this case.
Since $C^\ast(K)$ is a nuclear $C^\ast$\nbd{}algebra, it is an
amenable Banach algebra, so that $\HH^n(C^\ast(K))=0$ for $n\ge1$
(see~\cite{Khalkhali:Amenable}).  The space of bounded traces on
$C^\ast(K)$ is isomorphic to the space $\ell^1(\hat{K},\dim)$, where
$\dim$ denotes the measure that gives an $n$\nbd{}dimensional
irreducible representation volume~$n$.  Hence
$$
\HP^0(C^\ast(K)) \cong \HH^0(C^\ast(K)) \cong \ell^1(\hat{K},\dim),
\qquad
\HP^1(C^\ast(K))=\{0\}.
$$
We remark that Khalkhali shows that the entire and periodic cyclic
cohomology agree for amenable Banach algebras like $C^\ast(K)$.
However, he works with the von Neumann bornology.  The story is
totally different for the precompact bornology!

\section{Computations in analytic cyclic cohomology}
\label{sec:analytic}

In all relevant applications, we have $\HA_\ast \cong \HL_\ast$.  Such
an isomorphism should be expected because $\HO_\ast (\varinjlim C_i) =
\varinjlim \HO_\ast (C_i)$, that is, the homology of an inductive
system of complexes is already local.

However, the analytic and local theory also involve a completion
and the completion in the category of bornological vector spaces may
be badly behaved and, in particular, non-local.  The problem is that
the naïve completion need not be separated, so that we may have to
divide out a closure of $\{0\}$.  In general, this closure may be
surprisingly big and difficult to describe locally.  Nevertheless,
essentially the only completions we need are completions of tensor
products of spaces that are already complete.  Usually nothing goes
wrong with such completions.  Problems can only come from a very
serious failure of Grothendieck's approximation property.

If our bornological algebra is nuclear (or, more generally, satisfies
the bornological analogue of the approximation property), or if it is
a Fréchet algebra equipped with the precompact bornology, then the
completions are automatically separated, so that no problems appear.
We will explain this in a forthcoming article.

The algebra $\Rat(K)$ is nuclear, $\CINF(K)$ is both nuclear and
Fréchet and $C^\ast(K)$ is a Fréchet algebra.  As a result, we have
$\HA_\ast \cong \HL_\ast$ in these cases, so that
$$
\HA_\ast(\Rat(K)) \cong \HA_\ast (\CINF(K)) \cong \HA_\ast (C^\ast(K))
\cong \HL_\ast(C^\ast(K)).
$$
It is surprising that analytic and periodic cyclic homology are
different for $\CINF(K)$ despite the finite homological dimension.

Since ordinary cohomology is not local, there is no analogue of this
in cohomology.  In fact, the computation of analytic cyclic cohomology
is usually quite difficult.  For $\Rat(K)$, we get through with
Proposition~\ref{pro:sums} alone.  We see that the periodic, analytic
and local theories agree.  For the bigger algebras $\CINF(K)$ and
$C^\ast(K)$, we cannot compute the analytic cyclic cohomology.  We
only show that it is bigger than the periodic theory.  This is
straightforward using the Chern-Connes character constructed
in~\cite{Meyer:Analytic}.

It is well-known that in the entire theory there is a Chern-Connes
character for $\theta$\nbd{}summable Fredholm modules
(see~\cite{Connes:Entire}).  In fact, this was Connes's motivation to
introduce the theory.  In~\cite{Meyer:Analytic}, a different
construction is exhibited that produces a Chern-Connes character for
arbitrary Fredholm modules without summability conditions, that is, a
natural transformation
$$
\Chern \colon \K^\ast(A) \to \HA^\ast(A)
$$
for separable $C^\ast$\nbd{}algebras~$A$.  Its composition with the
natural map $\HA^\ast(A)\to\HL^\ast(A)$ is the Chern-Connes character
in local cyclic cohomology constructed by Puschnigg.

The Chern-Connes character can only exist if~$A$ is equipped with the
precompact bornology.  If~$A$ is an amenable Banach algebra, then the
result of Khalkhali mentioned above shows that there can be no
interesting Chern-Connes character with values in the analytic cyclic
cohomology of~$A$ equipped with the von Neumann bornology.

\begin{proposition}  \label{pro:HE_onto_HL}
  The canonical maps
  $$
  \HA^\ast(C^\ast(K))
  \to
  \HL^\ast(C^\ast(K))
  \cong
  \HL^\ast(\CINF(K))
  \leftarrow
  \HA^\ast(\CINF(K))
  $$
  are surjective.
\end{proposition}

\begin{proof}
  Since $\HL^1=0$ in both cases, we only have to consider the even
  cohomology.  The group $\K^0(A)$ is isomorphic to the space
  $\Hom(\Rep(K),\Z)$ by the universal coefficient theorem.  We
  identify this with the group $\Z^{\hat{K}} = \prod_{\rho\in\hat{K}} \Z$ of
  functions $\hat{K}\to\Z$.  Since analytic cyclic cohomology is a
  vector space over~$\C$, the Chern-Connes character yields a map
  $$
  \Chern \colon \Z^{\hat{K}} \otimes \C
  \cong
  \K^0(C^\ast(K)) \otimes \C
  \to
  \HA^0(C^\ast(K))
  \to
  \HA^0(\CINF(K)).
  $$
  If we compose these maps with the natural map to the local cyclic
  cohomology, we get the canonical map $\Z^{\hat{K}}\otimes\C \to
  \prod \C=\C^{\hat{K}}$.  It is easy to see that a function
  $\hat{K}\to\C$ belongs to its image if and only if its entries are
  contained in some finitely generated subgroup of~$\C$.

  Traces on $C^\ast(K)$ generate another subspace of
  $\HA^0(C^\ast(K))$, which is mapped onto the space
  $\ell^1(\hat{K},\dim)$.  Since~$\C$ has dense finitely generated
  subgroups, any sequence in~$\C$ can be decomposed as $(a_n) + (b_n)$
  with a sequence~$(a_n)$ whose entries lie in a finitely generated
  subgroup and a sequence $(b_n)\in\ell^1(\hat{K},\dim)$.  Therefore,
  the map $\HA^0(C^\ast(K)) \to \HL^0(C^\ast(K))$ is surjective.  This
  holds \emph{a fortiori} also for $\CINF(K)$.
\end{proof}

It is quite plausible that analytic and local cyclic cohomology should
agree for $C^\ast(K)$ or at least for $\CINF(K)$.  However, the author
does not know how to prove this.  Anyway, the proposition shows that
$\HA^\ast(\CINF(K))$ is bigger than $\HP^\ast(\CINF(K))$, that is, the
natural map $\HP^\ast(\CINF(K)) \to \HA^\ast(\CINF(K))$ is not
surjective.  This is striking because $\CINF(K)$ has finite
homological dimension.

\section{Concluding remarks}
\label{sec:conclusion}

If the analytic or local cyclic theories differ from the periodic
theory, then they yield results closer to $\K$\nbd{}theory.  This is
particularly obvious from the Universal Coefficient Theorem for the
local theory.  Hence it would be interesting in connection with the
Baum-Connes conjecture to compute the local or the analytic cyclic
homology for group algebras (both are isomorphic).

It is hard to conceive of a definition of finite dimension for
topological or bornological algebras that excludes the biprojective
algebra $\CINF(K)$.  Hence we should expect differences between the
periodic and analytic theory also for finite dimensional algebras, at
least if their $\K$\nbd{}theory is not finitely generated.  It seems
plausible that the existence of a finitely summable real spectral
triple with sufficiently good properties should guarantee equality
between the three cyclic theories.  If this spectral triple is a
``fundamental class,'' we get Poincaré duality between
$\K$\nbd{}theory and $\K$\nbd{}homology.  While this is a reasonable
assumption for non-commutative manifolds, it excludes algebras with
infinitely generated $\K$\nbd{}theory like the convolution algebras of
compact Lie groups studied above.  Hence such a concept of finite
dimensionality may not be sufficiently widely applicable.

\end{document}